\documentclass[11pt]{article}   
\usepackage{amssymb,amscd,latexsym}   
\usepackage{amsmath}
\usepackage{amsthm}
\newcommand{\rar}{\rightarrow}
\newcommand{\lar}{\longrightarrow}
\newcommand{\llar}{-\kern-5pt-\kern-5pt\longrightarrow}
\newcommand{\surjects}{\twoheadrightarrow}
\newcommand{\injects}{\hookrightarrow}

\newtheorem{Theorem}{Theorem}[section]
\newtheorem{Lemma}[Theorem]{Lemma}
\newtheorem{Corollary}[Theorem]{Corollary}
\newtheorem{Proposition}[Theorem]{Proposition}
\newtheorem{Remark}[Theorem]{Remark}
\newtheorem{Example}[Theorem]{Example}

\newtheorem{Definition}[Theorem]{Definition}
\newtheorem{Question}[Theorem]{Question}
\def\sqr#1#2{{\vcenter{\hrule height.#2pt
        \hbox{\vrule width.#2pt height#1pt \kern#1pt
            \vrule width.#2pt}
        \hrule height.#2pt}}}
\def\phi{\varphi}
\def\demo{\noindent{\bf Proof. }}
\def\square{\mathchoice\sqr64\sqr64\sqr{4}3\sqr{3}3}
\def\qed{\hspace*{\fill} $\square$}

\def\xx{{\bf x}}
\def\yy{{\bf y}}
\def\zz{{\bf z}}

\def\ff{{\bf f}}

\def\ff{{\bf f}}
\def\gg{{\bf g}}

\def\ii{\'{\i}}

\def\hht{{\rm ht}\,}

\def\Tor{\rm Tor}

\def\Reg{{\rm reg}}
\def\beg{{\rm indeg}}

\def\pp{{\mathbb P}}

\def\target{D\kern -0.001pt D}

\begin{document}

\begin{center}

{\Large{\bf\sc Implicitization of de Jonqui\`eres parametrizations}}
\footnotetext{Mathematics Subject Classification 2010
 (MSC2010).  13A30, 13H15, 13D02,13D45, 14E05,14E07.}

\vspace{0.3in}

{\large Seyed Hamid Hassanzadeh}\footnote{Partially supported by a  Post-Doc Fellowship (CNPq, Brazil).}
\quad\quad\
 {\large Aron  Simis}\footnote{Partially
supported by a CNPq grant.}

\end{center}


\begin{abstract}
One introduces a class of projective parameterizations that resemble generalized de Jonqui\`eres maps.
Any such parametrization defines a birational map $\mathfrak{F}$ of $\pp^n$ onto a hypersurface
$V(F)\subset \pp^{n+1}$ with a strong handle to implicitization. From this side, the theory here developed
extends recent work of Ben\ii tez--D'Andrea on monoid parameterizations.
The paper deals with both ideal theoretic and effective aspects of the problem.
The ring theoretic development gives information on the Castelnuovo--Mumford
regularity of the base ideal of $\mathfrak{F}$.
From the effective side, one gives an explicit
formula of $\deg(F)$ involving data from the inverse map of $\mathfrak{F}$ and show how the present parametrization
relates to monoid parameterizations.
\end{abstract}

\section{Introduction and notation}
\label{Pre}

Let $k$ denote an arbitrary infinite field which will be assumed to be algebraically closed for
the geometric purpose.
A rational map $\mathfrak{F}:\pp^n\dasharrow \pp^m$ is defined by $m+1$ forms $\ff=\{f_0,\ldots, f_m\}
\subset R:=k[\xx]=k[x_0,\ldots,x_n]$ of the same degree $d\geq 1$, not all null.
We often write $\mathfrak{F}=(f_0:\cdots :f_m)$ to underscore the projective setup.

The image of $\mathfrak{F}$ is the projective subvariety $W\subset \pp^m$ whose homogeneous
coordinate ring is the $k$-subalgebra $k[\ff]\subset R$ after degree renormalization.
Write $S:=k[\ff]\simeq k[\yy]/I(W)$, where $I(W)\subset k[\yy]=k[y_0,\ldots,y_m]$ is the homogeneous defining ideal
of the image in the embedding $W\subset \pp^m$.

We say that $\mathfrak{F}$ is {\em birational onto the image} if there is a rational map
backwards $\pp^m\dasharrow \pp^n$ such that the residue classes $\ff'=\{f'_0,\ldots, f'_n\}
\subset S$ of its defining coordinates do not simultaneously vanish and satisfy the
relations
\begin{equation}\label{birational_rule}
(\ff'_0(\ff):\cdots :\ff'_n(\ff))=(x_0:\cdots :x_n), \;
(\ff_0(\ff'):\cdots :\ff_m(\ff'))\equiv (y_0:\cdots :y_m)\mod I(W)
\end{equation}
Let $K$ denote the field of fractions of $S=k[\ff]$.
Note that the set of coordinates $(f'_0:\cdots :f'_n)$ defining the ``inverse'' map is not uniquely defined;
any other set $(f''_0:\cdots :f''_n)$ related to the first through requiring that it defines the same
element of the projective space $\pp^n_{K}=\pp^n_k\otimes_k {{\rm Spec}(K)}$ will do as well -- both tuples are called
{\em representatives} of the rational map (see \cite{bir2003} for details).
If $k$ is algebraically closed, these relations
translate into the usual geometric definition in terms of invertibility of the map on a dense Zariski open
set.

A special important case is that of a {\em Cremona map}, that is, a birational map
$$\mathfrak{G}=(g_0:\cdots :g_n):\pp^n\dasharrow\pp^n$$
of $\pp^n$ onto itself.
We assume, as usual, that the coordinates have no proper common factor.
In this setting, the common degree $d\geq 1$ of these coordinates is called the
{\em degree} of $\mathfrak{G}$.
Having information about the inverse map -- e.g., about its degree -- will be quite relevant in the sequel.
Thus, for instance, the structural equality
\begin{equation}\label{structure}
(g_0(g'_0,\ldots,g'_n):\cdots : g_n(g'_0,\ldots,g'_n))= (y_0:\cdots :y_n),
\end{equation}
involving the inverse map gives a uniquely defined form $D\in R$ such that $g_i(g'_0,\ldots,g'_n)=y_iD$, for every $i=1,\ldots,n$.
We call $D\in k[\yy]$  the {\em target inversion factor} of $\mathfrak{G}$.
By symmetry, there is a  {\em source inversion factor} $C\in k[\xx]$.

Our basic reference for the above is \cite{bir2003}, which contains enough of the introductory
material in the form we use here (see also \cite{AHA} for a more general overview).

\medskip

Now, the problem envisaged in this paper emerges from a particular situation of rational maps,
known as {\em elimination}.
Namely, one takes $m=n+1$ and assumes that $\dim k[\ff]=\dim R \,(=n+1)$.
Therefore, $W$ is a hypersurface defined by an irreducible form $F\in k[\yy]=k[y_0,\ldots,y_{n+1}]$.
We speak of $F$ informally as the {\em implicit equation} of $\mathfrak{F}$.
Elimination theory in this formulation is the problem of determining $F$ or at least its properties,
such as its degree. The set of the given forms defining $\mathfrak{F}$ is called
a {\em parametrization} of $F$. The theory has an applicable side shown in a very active research
area -- we refer to some of the
related modern work on the subject in the bibliography.

Although the main interest classically focused on {\em implicitization} -- i.e., in deriving the
implicit equation $F$ -- more recently quite some literature has appeared on the ideal theoretic structure
of the parametrization and the algebras naturally involved (\cite{CorDandr, BCS, Cox08, CHW, syl1, syl2}).
In this regard, a source of inspiration has been the classical {\em Sylvester forms},
a slightly imprecise notion to refer to certain generators
of the defining ideal of the Rees algebra associated to the base ideal of the rational
map $\mathfrak{F}$ (i.e., the ideal generated by the parameterizing forms).

Actually we go even more special, by dealing with rational maps which, in a sense, are allusive
of the classical de Jonqui\`eres plane Cremona map.
Namely, the class of parametrizations used here are suggestive of the stellar Cremona maps
by Pan (\cite{PanStellar}), a bona fide
generalization of the classical plane de Jonqui\`eres maps, and inspired by the results
of Ben\ii tez--D'Andrea (\cite{CorDandr}) on  the so-called {\em monoid
parametrizations}.

Precisely, start with a Cremona map $\mathfrak{G}=(g_0:\cdots :g_n):\pp^n\dasharrow\pp^n$
as explained above.
Let $f,g\in R$ be additional forms of arbitrary degrees  $\mathfrak{d}\geq 1$ and $d+ \mathfrak{d}$,
respectively.
We assume throughout that $f$ and $g$ are relatively prime.

\begin{Definition}\rm
The rational map $\mathfrak{F}=(g_0f:\cdots :g_nf:g):\pp^n\dasharrow\pp^{n+1}$ will be
called a {\em de Jonqui\`eres parametrization}.
\end{Definition}

Note the easy, though important, fact that $\mathfrak{F}$ is a birational map onto its
image $W=V(F)$.
This follows immediately from the usual field extension criterion (see, e.g., \cite[Proposition 1.11]{AHA}.
Moreover, if the inverse of $\mathfrak{G}$ is $\mathfrak{G}^{-1}=(g_0':\cdots :g_n')$,
with $g_i'\in k[y_0,\ldots,y_n]$, then $(\overline{g_0'}:\cdots :\overline{g_n'})$ is a representative
 of $\mathfrak{F}$, where the bar over an element of $k[\yy]$
denotes its class modulo $(F)$ -- note that this representative of $\mathfrak{F}^{-1}$
does not involve the last variable $y_{n+1}$.

The Cremona map $\mathfrak{G}$ may be called the {\em underlying} (or {\em structural}) Cremona map of  $\mathfrak{F}$.

\medskip

The main results of the paper are stated in Theorem~\ref{pregularity}, Proposition~\ref{pdegree},
Proposition~\ref{Rees_equations} and Theorem~\ref{From_monoid}.

Let us briefly describe the contents of the next sections.

Section~\ref{first} gives the main properties of the base ideal of the parametrization, such as
structure of syzygies, free resolution and regularity.
Part of the information of this section is crucial for introducing the concept of {\em syzygetic polynomials}
that arise as natural candidates for the implicit equation (often with extraneous factors).

Section~\ref{second} deals with  the implicit equation $F$. Here one introduces the basic
polynomials that play a role in the nature of $F$, such as the syzygetic polynomials mentioned before.
One heavily draws on the hypothesis that the de Jonqui\`eres parametrization is birational, by having
the defining parametrization of the inverse map and the inversion factor take control of the situation.
The section also examines the details of two main cases of the given de Jonqui\`eres parametrization, called
``the inclusion case'' and the ``non-zero-divisor case'', respectively.
It is worth pointing that the first of these two cases covers as a very special case the situation of
a monoid parametrization.

In Section~\ref{third} one focus on the so-called ``Rees equations'' of the parametrization.
These are the elements of a minimal set of generators of a presentation ideal (the ``Rees ideal'')
of the Rees algebra of the base ideal of
$\mathfrak{F}$, one of which, of course, is $F$ itself.
These have been variously studied by several authors, some listed in the references.
The idea in this section is based on the method of downgrading that has been used in different sources
(e.g., \cite{BCS}, \cite{Trento}, \cite{syl2}).
Ours is a modification of this method -- hereby called {\em birational downgrading} -- by which we use
the forms defining the inverse map rather than the usual procedures in the literature.
The main result yields a set of Rees equations candidates for a set of minimal generators,
generating an ideal having as a minimal prime component the entire Rees ideal.
The sections end with a result giving the precise relation between the Rees ideal of de Jonqui\`eres
parameterizations and the one of the monoid parameterizations.

\section{Syzygetic background}\label{first}

In this section we establish the basic relations of degree $1$ of the forms $g_0f,\ldots,g_nf, g$
defining the rational map. For the next lemma and proposition, $(g_0:\cdots :g_n)$ defines
any rational map, not necessarily Cremona.

\subsection{A mapping cone}\label{on_mapping_cone}

In this part we state a very general result regarding a certain mapping cone naturally associated to the
present data.
The construction is completely general and does not require a graded situation.
Accordingly, we refresh our data just assuming that $I\subset R:=k[x_0,\ldots,x_n]$ is an arbitrary
ideal and $f,g\in R$ are given elements.

\begin{Lemma}\label{basic_transfer} If $\gcd(f,g)=1$ then:
\begin{enumerate}
\item[{\rm (a)}] $If:(g)=(I:(g))f$.
\item[{\rm (b)}] Multiplication by $g$ induces an isomorphism $R/(I:(g))f\simeq (If,g)/If$ of $R$-modules.
\end{enumerate}
\end{Lemma}
\demo (a) The inclusion $If:(g)\supset (I:(g))f$ is obvious regardless of any relative assumption about $f,g$.
Conversely, let $b\in R$ be such that $bg\in If$. Then $f$ divides $bg$ and, since $\gcd(f,g)=1$, then
$f$ divides $b$. Say, $b=af$, with $a\in R$. Then $(ag)f\in If$, hence $ag\in I$, i.e., $a\in I:(g)$.
Therefore, $b\in (I:(g))f$.

\smallskip

(b) One has $(If,g)/If\simeq (g)/(g)\cap If=(g)/(If:(g))g\simeq R/If:(g)$, where the last isomorphism is
multiplication by $g^{-1}$.
Now apply (a).
\qed

\medskip

Quite generally, a surjective $R$-module homomorphism $\pi:R^q\surjects I:(g)$ induces a {\em content map}
$c(g):R^q\lar R^{n+1}$. In explicit coordinates: let $\pi$ be induced by choosing a set of generators $\{c_1,\ldots,c_q\}$
of $I:(g)$, so that $\pi(v_j)=c_j$, where $\{v_1,\ldots,v_q\}$ is the canonical basis of $R^q$.
Given a set $\{g_0,\ldots,g_p\}$ of generators of $I$, let
$\{e_0,\ldots,e_p\}$ denote the canonical basis of $R^{p+1}$.
Write $c_jg=\sum_{i^=0}^p h_{ij}g_i$, with $h_{ij}\in R$.
Then $c(g)(v_j)=\sum_{i^=0}^p h_{ij}e_i$, for $j=1,\ldots,q$.

This simple construction will be used in the following result.

\begin{Lemma}\label{mapping_cone}
Let $\mathfrak{R}$ and $\mathfrak{S}$ denote finite free resolutions of $R/I$ and $R/(I:(g))f$,
respectively.
Then multiplication by $g$ lifts to a map $\mathfrak{S}\rar \mathfrak{R}$ whose associated mapping
cone is a free resolution of $R/(If,g)$.
In particular, a syzygy matrix of the generators of $(If,g)$ has the form
$$\Psi= \left(
\begin{array}{cc}
\phi & c(g)\\
\mathbf{0} & -f\pi
\end{array}
\right),
$$
where $\phi$ denotes a syzygy matrix of a given set of generators of $I$.
\end{Lemma}
\demo
As in Lemma~\ref{basic_transfer}(b), multiplication by $g$ induces an injective $R$-module homomorphism
$R/(I:(g))f\injects R/If$ with image $(If,g)/If$.
This homomorphism lifts to a map of complexes (free resolutions)
$$
\begin{array}{cccccccccccc}
\mathfrak{R}: & \cdots & \lar & R^{m_1} &\stackrel{\phi}{\lar} & R^{p+1} & \stackrel{f\gg}{\lar} & R & \lar & R/If &\rar & 0\\[6pt]
&& &             \uparrow   &                 & \kern-12ptc(g)\uparrow         &         &\kern-8pt\cdot g\uparrow  && \cdot g\uparrow &&\\[6pt]
\mathfrak{R}: &  \cdots  & \lar & R^{r_1} &\stackrel{\psi}{\lar}              & R^q     & \stackrel{f\pi}{\lar} & R & \lar & R/(I:(g))f &\rar & 0
\end{array},
$$
where $\gg=(g_0\cdots g_p)$, with  and $c(g)$ is the above content map.
Then the corresponding mapping cone is an $R$-free resolution of $(R/If)/((If,g)/If)\simeq R/(If,g)$ (see \cite[Exercise A3.30]{E}).
\qed

\subsection{Graded minimality and regularity}

We move back to the original graded situation.
Namely, set $I=(g_0,\ldots,g_n)$, where the $g_i$'s are forms of degree $d\geq 1$ minimally generating $I$, and $f,g$ are forms
with $\deg(g)=d+\deg(f)$ such that $\gcd(f,g)=1$.
Also, let $\{c_1,\ldots,c_q\}$ be a set of minimal generators of $I:(g)$, with $c_j$ homogeneous of degree $C_j$.

Let
$$\cdots\rar\bigoplus_{j=1}^{m_1} R(-a_{1j})\xrightarrow{\phi} \bigoplus_{i=0}^{n} R(-d) \xrightarrow{\gg}R \rightarrow R/I\rightarrow 0$$
and
 $$\cdots\rar\bigoplus_{j=1}^{q_2} R(-C_{2j})\xrightarrow{\psi}\bigoplus_{j=1}^{q} R(-C_{j})\xrightarrow{\pi}R \rightarrow R/I:g\rightarrow 0$$
 stand for minimal graded free resolutions of $R/I$ and $R/I:g$, respectively,
from which we immediately derive minimal graded free resolutions of $R/If$ and $R/(I:g)f$:
$$\cdots \bigoplus_{j=1}^{m_1} R(-a_{1j}-\deg(f))\xrightarrow{\phi_1=\phi} R(-(d+\deg(f)))^n
\xrightarrow{f\,\gg}R \rightarrow R/If\rightarrow 0 ,$$
 $$\cdots\rar\bigoplus_{j=1}^{q_1} R(-C_{1j}-\deg(f))\xrightarrow{\psi_1=\psi}\bigoplus_{j=1}^{q}
 R(-C_{j}-\deg(f))\xrightarrow{f\,\pi}R \rightarrow R/(I:g)f\rightarrow 0.$$
Shifting the second of these resolutions by $-(d+\deg(f))$, one obtains a map of complexes, where the vertical
 homomorphisms are also homogeneous of degree $0$

 {\scriptsize
 $$
 \begin{array}{ccccccccccccc}
 \cdots\kern-6pt &\rar \kern-6pt & \bigoplus_{j=1}^{m_{i}} R(-a_{ij}-\mathfrak{d})\kern-6pt &\rar
 &\kern-6pt \cdots \kern-6pt & \stackrel{\phi_2}{\rar}\kern-6pt  &  R(-(d+\mathfrak{d}))^n
 \kern-6pt & \rar\kern-6pt & \kern-6pt R \kern-6pt & \rar   \kern-6pt& R/If  \kern-6pt &\rar  & 0 \\
  && c_i(g)\uparrow &&&& c(g)\uparrow && \cdot g\uparrow && \cdot g\uparrow &&\\
   \cdots\kern-6pt &\rar \kern-6pt & \bigoplus_{j=1}^{q_{i}} R(-C_{ij}-(d+2\mathfrak{d}))\kern-6pt &\rar
   &\kern-6pt \cdots \kern-6pt & \kern-6pt\stackrel{\psi_2}{\rar}  &
 \kern-6pt \bigoplus_{j=1}^{q} R(-C_{j}-(d+2\mathfrak{d})) &\kern-6pt \rar & \kern-6pt R(-(d+\mathfrak{d}))
 & \rar \kern-6pt  & \frac{R}{(I:g)f}(-(d+\mathfrak{d})) & \kern-6pt\rar  & 0
 \end{array}
 $$
 }
 where we have written $\mathfrak{d}:=\deg(f)$ for editing purpose.

 \smallskip

We let $\Reg(M)$ denote the Castelnuovo--Mumford regularity of a graded $R$-module and let hd$(M)$ stand for its homological
(i.e., projective) dimension.

\smallskip

 One has:
 \begin{Proposition}\label{mapping_cone_graded}
 With the above notation, the associated mapping cone $\mathcal{C}_{\bullet}$ is
  a graded free resolution of $R/(If,g)$:
  {\scriptsize
 \begin{equation}\label{resolution_of_stellar}\nonumber
\cdots\rar \left(\bigoplus_{j=1}^{m_1} R(-a_{1j}-\deg(f))\right)\oplus
\left(\bigoplus_{j=1}^qR(-C_j-(d+2\deg(f)))\right)\xrightarrow{\Psi}
R(-(d+\deg(f)))^{n+1}\rightarrow R \rightarrow  R/(If,g)\rightarrow 0.
\end{equation}
}
Moreover, if $\Reg(R/I)\leq d+\deg(f)-2$ then
this resolution is minimal.
\end{Proposition}
\demo
Applying the minimality criterion stated in \cite[Exercise A3.30]{E} it suffices to show that
$-a_{ij}-\deg(f) >-C_{ik}-(d+2\deg(f))$ for all $i, j,k$. Now, on one hand,
$a_{ij}-i$ is at most the regularity of $R/I$, for any $i,j$; on the other hand,
for any $i,k$, $C_{ik} \geq i-1$, where $C_{0k}=C_k$. Therefore the condition is fulfilled if the
regularity of $R/I$ is bounded as stated.
\qed

\begin{Corollary}\label{hom_dim}
With the above notation, assume that  hd$(R/(I:g)f)\leq {\rm hd}(R/I)-1$
{\rm(}e.g., if $g\in I$ and $I$ has codimension $\geq 2${\rm)}.
Then
$${\rm hd}(R/(If,g))\leq {\rm hd}(R/I),$$
with equality provided $\Reg(R/I)\leq d+\deg(f)-2$.
\end{Corollary}

\medskip

Since the preceding regularity bound implies, in particular, the minimality of the above graded free presentation
of $R/(If,g)$, thus having a direct impact on the search for a minimal set of bihomogeneous Rees equations,
it is pertinent to understand how this bound reflects on the current data.

\begin{Proposition}\label{reg_of_stellar} Keeping the previous notation, one has:
\begin{itemize}
\item[{\rm (a)}] $\Reg(R/(If,g))\leq \max\{reg(R/I)+\deg(f),\Reg(R/(I:g))+d+2\deg(f)-1\}$.
\item[{\rm (b)}] If\,, moreover, $\Reg(R/I)\leq d+\deg(f)-2$ then
$$\Reg(R/(If,g))= \Reg(R/(I:g))+d+2\deg(f)-1.$$
\end{itemize}
\end{Proposition}
\demo
 (a)  Computing the regularity in
 terms of the twists of the graded free resolution $\mathcal{C}_{\bullet}$ in {Proposition}~\ref{mapping_cone_graded},
 one finds
\begin{eqnarray*}\Reg(R/(If,g))&\leq & \max\{\Reg(R/If),\Reg(R/(I:g)f)+ d+\deg(f)-1\}\\
&= &\max\{\Reg(R/I)+\deg(f),\Reg(R/(I:g))+d+2\deg(f)-1\}.
\end{eqnarray*}

 (b) By the second assertion in {Proposition}~\ref{mapping_cone_graded},
 the mapping cone is a graded minimal free resolution. Therefore,  if $\Reg(R/I)\leq d+\deg(f)-2$ then
 the maximum in (a) is  the second term.
\qed

\subsection{Regularity in the case of isolated base points}

 We keep the notation of the previous subsection.
 Namely, $I=(g_0,\ldots,g_n)$, where the $g_i$'s are forms of degree $d\geq 1$ minimally generating $I$, and $f,g$
 are nonzero forms such that $\deg(g)=d+\deg(f)$ and $\gcd(f,g)=1$.
 For any ideal $\mathfrak{a}\subset R$, we denote by $\mathfrak{a}^{\rm sat}$ its saturation $\mathfrak{a}:(\xx)^{\infty}$.
 If $M$ is a graded $R$-module 
we will set
$$\beg (M):=\inf \{ \mu \ \vert \ M_\mu \not= 0\},$$
with the convention that $\beg (0)=+\infty$, and
$${\rm end} (M):= \sup \{ \mu \ \vert \ M_\mu \not= 0 \},$$ 
with the  convention that ${\rm end}(0)=-\infty$.

\begin{Theorem}\label{pregularity} Suppose that $\dim(R/I)\leq 1$. Then
\begin{itemize}
\item[{\rm (1)}] $\Reg(R/I)=\max\{(n+1)(d-1)-\beg(I^{\rm sat}/I),n(d-1)-\beg((\alpha):I/I)\}$,
where $\alpha$ denotes a maximal regular sequence of $d$-forms in $I$.
\item[{\rm (2)}] If in addition $I$ is the base ideal of a Cremona map, then $\Reg(R/I)\leq n(d-1)-1.$
\item[{\rm (3)}] $\Reg(R/(If,g))\leq \Reg(R/I)+\deg(f)+\deg(g)-1$, with equality holding
when $g$ is a non-zero-divisor on $R/I$.
\item[{\rm (4)}] If  $\Reg(R/I)\leq \deg(g)-2$ then  $\Reg(R/I:g)\leq \Reg(R/I)$.
\end{itemize}
\end{Theorem}
\demo  (1) We copy {\em ipsis litteris} the argument in the proof of \cite[Theorem 1.5]{HS}, updating the setup.
Thus, the ground ring has now dimension $n+1$, the base ideal $I\subset R$ has codimension $n$ (hence,
 $I^{\rm un}=I^{\rm sat}$ and, necessarily, $d\geq 2$).
 Note that one can always pick a maximal regular sequence of $d$-forms in $I$.
 We thus obtain $(I^{\rm sat}/I)^{\Check{}}=(I^{\rm sat}/I)((n+1)(d-1))$, hence ${\rm end}(I^{\rm sat}/I)+\beg(I^{\rm sat}/I)=(n+1)(d-1)$.

\medskip

(2) In the case where $I$ is the base ideal of a Cremona map,  $\beg(I^{sat}/I)\geq d+1$ according to \cite{PanRusso} (by convention,
if $I$ is saturated one sets $\beg(I^{sat}/I)=+\infty$). Therefore $(n+1)(d-1)-\beg(I^{sat}/I)\leq n(d-1)-2$. On the other hand $(\alpha)\subsetneq I$
since $I$ defines a Cremona map and $d\geq 2$, where $\alpha$ is an in item (1).
This gives  $\beg((\alpha):I/I)\geq 1$ which implies that   $n(d-1)-\beg((\alpha):I/I)\leq n(d-1)-1$. The assertion then follows from (1).

\medskip

 (3)
 Let us argue that $\dim({\Tor}_i^R(R/If,R/g))\leq 1$ for all $i$. In the present situation, the relevant indices are
 $i=0,1$. Now,  ${\Tor}_0^R(R/If,R/g))=R/(If,g)$ and
 \begin{eqnarray*}
 {\Tor}_1^R(R/If,R/g))&\simeq & If\cap (g)/(g) If\simeq (If:g)g/If\cdot g
 \simeq  (If:g)/If(-\deg(g))\\ &\simeq &(I:g)/I(-\deg(f)-\deg(g)),
 \end{eqnarray*}
 whose annihilator is the ideal $I:(I:g)$.
 Thus, both modules have dimension at most one as $\dim(R/I)\leq 1$ by assumption.
Therefore, we can apply \cite[Theorem 0.2 or Corollary 5.8]{C} which gives

\begin{eqnarray}\label{torreg}
&& \max \{\Reg(R/(If,g)), \Reg(I:g/I)+\deg(f)+\deg(g)-1\}\\\nonumber
&=&\max\{\Reg(R/(If,g)), \Reg((If:g)/If)-1\}\\\nonumber
&=&\max\{{\Tor}_0^R(R/If,R/g), {\Tor}_1^R(R/If,R/g)-1\}
=\Reg(R/I)+\deg(f)+\deg(g)-1.
\end{eqnarray}
Obviously, the value in the first line of (\ref{torreg}) bounds $\Reg(R/(If,g))$, hence we are through.
If $g$ is a non-zero divisor on $R/I$ then $I:g/I=0$, hence the maximum value in the first line of (\ref{torreg})
is $\Reg(R/(If,g))$.

\medskip

(4) This follows from item (3) and Proposition~\ref{reg_of_stellar} (b).
\qed

\section{The search for the implicit equation}\label{second}

We keep the notation of Section~\ref{Pre}.

\subsection{Monoids and syzygetic polynomials}

\medskip

A form $F$ of the  shape $F=G+Hy_{n+1}$, where $G,H$ are forms in $k[y_0,\ldots,y_n]$, is called a  {\em monoid}
(cf. \cite{Piene_et_al} for generalities on these forms).
Thus, this is simply a polynomial of degree $1$ in the one variable polynomial ring $B[y_{n+1}]$,
with homogeneous coefficients in $B=k[y_0,\ldots,y_n]$.
As we will see, a good deal of the results involves monoids throughout.
We will often say an $y_{n+1}$-monoid to stress the privileged variable $y_{n+1}$.

The following gadget will be basic throughout.
Consider a syzygy of $J=(If,g)$ as in Lemma~\ref{mapping_cone} with nonzero last coordinate.
Its polynomial version is a $1$-form in $R[y_0,\ldots,y_n,y_{n+1}]$
\begin{equation}\label{pol_syzygy}
\sum_{i^=0}^n h_{ij}y_i-fc_jy_{n+1},
\end{equation}
where $I:(g)=(\ldots,c_j,\ldots)$ and $c_jg=\sum_{i^=0}^n h_{ij}g_i$.

\begin{Definition}\label{syzygetic_pol}\rm Suppose that $(g_0:\cdots :g_n)$ defines
a Cremona map and let $(g'_0:\cdots :g'_n)$ define its inverse map.
The $j$th {\em syzygetic polynomial} is the form
\begin{equation}\label{syzygetic_pol_eq}
\sum_{i^=0}^n h_{ij}(g'_0,\ldots,g'_n)y_i-f(g'_0,\ldots,g'_n)c_j(g'_0,\ldots,g'_n)y_{n+1}\in
k[\yy]:=k[y_0,\ldots,y_n,y_{n+1}],
\end{equation}
obtained from (\ref{pol_syzygy}) by evaluating $x_i\mapsto g'_i$, for $i=0,\ldots,n$.
\end{Definition}

The main property of such forms is that they belong to the defining ideal of the Rees algebra of the ideal $(If,g)$.
To see this, recall that since the rational map defined by the generators of $(If,g)$ is birational onto $V(F)$,
one has a $k$-algebra isomorphism of the Rees algebras
\begin{equation}\label{Rees_are_isomorphic}
\mathcal{R}_R((If,g))\simeq \mathcal{R}_{S}(I'),
\end{equation}
where $I'=(\overline{g_0'},\ldots ,\overline{g_n'})$ and $S=k[\yy]/(F)$ (see
\cite[Theorem 2.18, proof of (a) $\Rightarrow$ (b)]{AHA}).
This isomorphism is induced by the identity map of $R[\yy]=k[\yy][\xx]=k[\xx,\yy]$.
In terms of the respective defining ideals  $\mathcal{J}$ and  $\mathcal{K}$ over $k[\xx,\yy]$, we
have an equality $\mathcal{J}=\mathcal{K}$.
Therefore, by definition of $\mathcal{K}$, the syzygetic polynomial
$$P=P(\yy):=\sum_{i=0}^n h_i(g'_0,\ldots,g'_n)y_i-f(g'_0,\ldots,g'_n)c_j(g'_0,\ldots,g'_n)y_{n+1}$$
vanishes modulo $Ik[\xx,\yy]$.
But since it is a polynomial in $\yy$ only,  it necessarily belongs to $(F)\subset k[\yy]$, i.e., it is a
multiple of the implicit equation $F$.

This suggests that a syzygetic polynomials is a fair candidate for the implicit equation $F$ and will coincide with $F$ up to
a nonzero field element provided it be irreducible.

\subsection{The degree of the implicit equation}

\medskip

In this part we establish a formula for the degree of $F$ in terms of the data introduced so far.

\begin{Proposition}\label{pdegree}
Suppose that $\mathfrak{G}=(g_0:\cdots :g_n):\pp^n\dasharrow\pp^n$
defines a Cremona map of degree $d$ with inverse map $\mathfrak{G}^{-1}=(g'_0:\cdots :g'_n)$ and
target inversion factor $D\subset k[y_0,\ldots,y_n]$.
Let $f,g\in R=k[x_0,\ldots,x_n]$ be forms, with $\deg(g)=d+\deg(f)$, such that
$f(\gg'), g(\gg')$ are nonzero as elements of $k[y_0,\ldots,y_n]$.
Letting $F\subset k[y_0,\ldots, y_n,y_{n+1}]$ denote the implicit equation of the parametrization
$(fg_0:\cdots :fg_n:g):\pp^n\dasharrow \pp^{n+1}$, one has:
\begin{enumerate}
\item[{\rm (i)}] $F$ is the $y_{n+1}$-monoid
\begin{equation}\label{elimination_candidate}
\frac{g(\gg')-y_{n+1}f(\gg')\,D}{\gcd(g(\gg'), f(\gg')\,D)}\,.
\end{equation}
\item[{\rm (ii)}] $\deg(F)=\deg(g)\deg(\mathfrak{G}^{-1})-\deg(\gcd(g(\gg'), f(\gg')\,D))=
\deg(f)\deg(\mathfrak{G}^{-1})+\deg(D)+1 -\deg(\gcd(g(\gg'), f(\gg')\,D))$.
\end{enumerate}
\noindent In particular, $\deg(F)\leq\deg(g)\deg(\mathfrak{G}^{-1})$.
\end{Proposition}
\demo Set $S:=k[y_0,\cdots,y_{n+1}]/(F)$ for the homogeneous coordinate ring of the image of $\mathfrak{F}$.
Since the associated de Jonqui\`eres parametrization defined by the generators of $(If,g)$
is birational, a formula such as (\ref{structure}) implies the vanishing of the $2\times2$ minors of
 $$
\begin{pmatrix}(fg_{_0})(\gg')&\cdots&(fg_{_n})(\gg')&g(\gg')\\
y_0&\cdots&y_n&y_{n+1}
\end{pmatrix}
$$
modulo $(F)$. In particular, each of the minors $P_i:=y_ig(\gg')-y_{n+1}f(\gg')g_i(\gg')$ fixing
the last column, are multiples of $F$.
On the other hand, using (\ref{structure}) for the Cremona map $\mathfrak{G}$ yields
 $g_i(\gg')=y_i\,D$, for $i=0,\ldots,n$.
Since $y_i$ is not a factor of $F$, it follows that $F$ divides the nonzero $y_{n+1}$-monoid $g(\gg')-y_{n+1}f(\gg')\,D$.
Clearly, $F$ is not a factor of $\gcd(g(\gg'), f(\gg')\,D)$ as the latter lives in $k[y_0,\ldots,y_n]$
while $F$ involves effectively the variable $y_{n+1}$.

Now
\begin{equation*}
\frac{g(\gg')-y_{n+1}f(\gg')\,D}{\gcd(g(\gg'), f(\gg')\,D)}
\end{equation*}
is an $y_{n+1}$-monoid with relatively prime components, hence is irreducible.
Therefore, it must coincide with $F$.

(ii) Taking degrees, the first equality in the stated formula follows readily, while
the subsequent equality follows from the standing assumption that $\deg(g)=\deg(f)+\deg(\mathfrak{G})$
and from the equality $\deg(\mathfrak{G})\deg(\mathcal{G}^{-1})=\deg(D)+1$ by definition of the
inversion factor $D$.
\qed

\begin{Corollary}\label{evaluated_gcd} If $\gcd(f(\gg '),g(\gg'))=1$ then
$$\deg(F)= \deg(g)\deg(\mathfrak{G}^{-1})-\deg(\gcd(g(\gg'),D)). $$
In particular,  $\deg(f)\deg(\mathfrak{G}^{-1})+1\leq \deg(F)< \deg(g)\deg(\mathfrak{G}^{-1})$.
\end{Corollary}
\demo
The displayed equality follows immediately from the formula in Proposition~\ref{pdegree},
so only the lower bound is the question.
For that,  writing $\deg(g)=\deg(f)+\deg (\mathfrak{G})$ yields
$\deg(F)= \deg(f)\deg(\mathfrak{G}^{-1})+\deg (\mathfrak{G})\deg (\mathfrak{G}^{-1})-\deg(\gcd(g(\gg'),D))$.
But $\deg (\mathfrak{G})\deg (\mathfrak{G}^{-1})=\deg(D)+1$, while obviously $\deg(\gcd(g(\gg'),D))\leq \deg(D)$.
\qed

\subsection{The inclusion case}

\medskip

We focus on the case where $g\in I$.

Here one has $I:(g)=R$, which in the notation
of Subsection~\ref{on_mapping_cone}
 tells us that $\pi:R\lar R$ can be taken to be the identity map
and the content map $c(g):R\lar R^{n+1}$ picks up only one additional syzygy.
Here the syzygy matrix of the given generators of $(If,g)$ is of the form
\begin{equation}\label{inclusion_matrix}
\Psi= \left(
\begin{array}{cc}
\phi & c(g)\\
\mathbf{0} & -f
\end{array}
\right),
\end{equation}
where $\phi$ denotes a syzygy matrix of the given set of generators of $I$
and $c(g)$ stands the column vector defining the content map.

 Thus, there is only one syzygetic polynomial $P:=\sum_{i=0}^n h_i(\gg')y_i-f(\gg')y_{n+1}$.
 Keeping the assumptions of {\rm Proposition~\ref{pdegree}}, one has:

\begin{Proposition}\label{inclusion}
Assume that $\gcd(f(\gg '),g(\gg'))=1$.
If $P\in k[y_0,\ldots,y_n,y_{n+1}]$ is a syzygetic polynomial, the following conditions are equivalent:
\begin{itemize}
\item[{\rm (i)}] $g\in I$.
\item[{\rm (ii)}] $\deg(F)= \deg(f)\deg(\mathfrak{G}^{-1})+1$ and  $(F)=(P)$.
\end{itemize}
\end{Proposition}
\demo (i) $\Rightarrow$ (ii)
Quite generally, when $g\in I$ one has $\deg(P)=\deg(f)\deg(\mathfrak{G}^{-1})+1$.
On the other hand, by Corollary~\ref{evaluated_gcd}, $\deg(F)\geq \deg(f)\deg(\mathfrak{G}^{-1})+1$.
Since $F$ is a factor of $P$, we are through.

(ii) $\Rightarrow$ (i)
Since $\deg(P)= \deg(f)\deg(\mathfrak{G}^{-1})+1$, confronting with the general shape of $P$
as in Definition~\ref{syzygetic_pol} yields $\deg (c_j(g'_0,\ldots,g'_n))=0$ for a generator $c_j$ of the
conductor $I:g$. This forces $c_j$ to be invertible, so $g\in I$.
\qed

A special case of Proposition~\ref{inclusion} is a result of \cite{CorDandr}:

\begin{Corollary}\label{Benitez_Dandrea}
If $\mathfrak{G}$ is the identity map of $\pp^n$ and $\gcd(f,g)=1$ then
\begin{equation}F= g(y_0,\ldots,y_n)-f(y_0,\ldots,y_n)\,y_{n+1}.
\end{equation}
In particular, $\deg(F)=\deg(f)+1$.
\end{Corollary}
\demo In this case, the inverse is also the identity, so $\gcd(f(\gg'),g(\gg'))=\gcd(f(\yy),g(\yy))=1$.
On the other hand, the target inversion factor is $1$, so the polynomial (\ref{elimination_candidate})
is $g(y_0,\ldots,y_n)-f(y_0,\ldots,y_n)\,y_{n+1}$.
\qed

\medskip

At the other end of the spectrum, so to say, we find as a consequence a more ``typical'' situation:

\begin{Corollary}\label{general_form}
Keeping the assumptions of $\,${\rm Proposition~\ref{pdegree}}, assume that $g\in I$.
If $f$ is a general form then $(F)=(P)$.
\end{Corollary}
\demo Since $f$ is chosen to be general and $g$ are $\gg'$ are fixed once for all,
 the condition $\gcd(f(\gg'),g(\gg'))=1$ is fulfilled.
 \qed

\medskip

If $f$ is not sufficiently general it may happen that $P$ as above is not irreducible as the following
example entails.

\begin{Example}\rm
Take the maximal minors of the following $4\times 3$ matrix over $R=k[x_0,x_1,x_2,x_3]$
$$\left(
\begin{array}{ccc}
0   &  0  &  -x_1\\
-x_0& x_0-x_1& x_1\\
x_0  & 0  &  0\\
x_2 & -x_3 &  x_3
\end{array}
\right).
$$
These $3$-forms $g_0,g_1,g_2,g_3$ define a Cremona map of $\pp^3$ with inverse given by the $2$-forms
$$-y_0y_3,\, y_0y_2,\, -y_1y_3-y_2y_3,\, -y_2^2-y_2y_3
$$
(\cite[Section 2.1]{Russimis}).
If one takes $g\in I=(g_0,g_1,g_2,g_3)$ and $f$ sufficiently special, but still such that
gcd$(f,g)=1$, then it is apparent that $y_3$ will
come out as a factor of $P$ -- e.g., take $g=x_0g_3=x_0^3x_3$ (the minor corresponding to the last $3$ rows) and $f=x_0+x_2$;
then $P=-y_3F$, where $F$ is the implicit quadric equation.
\end{Example}

Yet another special notable case of $g\in I$  is worth isolating as well, where the data are somewhat twisted around.
We recall that a form $g\in R=k[\xx]$ is called {\em homaloidal} if its partial derivatives
(the so-called {\em polar map} of $G$) define a Cremona map.
The ideal generated by the partial derivatives of a form is often called its {\em gradient ideal}.

\begin{Corollary}\label{Eulerian_elimination} {\rm (char$(k)=0$)}
Let $g\in R=k[\xx]$ denote a reduced homaloidal form of degree $d+1$ and let
$I\subset R$ stand for the gradient ideal of $g$, and let $\{g'_0,\ldots,g'_n\}\subset k[y_0,\ldots,y_n]$
define the inverse map of the polar map of $g$.
If $f=\sum_{i=0}^n\lambda_ix_i$ is a general linear form then
\begin{equation}\label{eulerian_F}
F= \sum_{i=0}^n(y_i-(d+1)\lambda_i\,y_{n+1})\,g'_i(y_0,\ldots,y_n).
\end{equation}
In particular, $\deg(F)= \deg(g'_i) +1$.
\end{Corollary}

The polynomial (\ref{eulerian_F}) might be called the {\em general Eulerian equation} of a polar Cremona map.

\begin{Remark}\rm
We note that under the hypothesis that $g\in I$, there is an inclusion $J\subset I$.
This triggers a natural injection $\mathcal{R}(J)\subset \mathcal{R}(I)$ of Rees algebras.
Thus, in principle, this would give information about the defining Rees equations of $J$
out of these of the base ideal $I$.
However, setting up explicit presentations requires moving around variables, so the ultimate
computational advantage is not so clear.
Also note that if, moreover, $I$ is saturated
and $f$ is sufficiently general, then $J=I\cap (f,g)$ and $J:I=(f,g)$ (to see the last equality,
note it is obvious if $\hht I\geq 3$ since $\{f,g\}$ is a regular sequence, and if $\hht I=2$ we just
need that no minimal prime of $I$ be a minimal prime of $(f,g)$, which is the case if $f$ is
sufficiently general).
One may ask how implicitization may profit from this simple situation of linkage in a {\em coarse} sense.
\end{Remark}

\subsection{The non-zero-divisor case}

\smallskip

Assume that $g$ is a non-zero-divisor on $R/I$.
In this situation, $I:(g)=I$, hence the map $\pi$ in Lemma~\ref{mapping_cone} boils down to the structural
surjection $\phi: R^{n+1}\lar I$.
Accordingly, the content map $c(g)$ reduces to $g$ times the identity map of $R^{n+1}$.
Therefore, a presentation matrix of $J=(If,g)$ has now the form
$$\Psi= \left(
\begin{array}{cc}
\phi & \kern4ptg\cdot \mathbf{1}_{n+1}\\
\mathbf{0} & -f\gg
\end{array}
\right),
$$
where $\phi$ is a syzygy matrix of $I$ and $\gg = (g_0\cdots g_n)$.

\begin{Proposition}\label{elimination_nzd_case}
Let $\mathfrak{G}=(g_0:\cdots :g_n):\pp^n\dasharrow\pp^n$
stand for a Cremona map of degree $d$ with base ideal $I=(g_0,\ldots,g_n)$
and let $f,g\in R$ be given as before.
Suppose that $g$ is a non-zero-divisor on $R/I$.
Then the implicit equation $F$  is a factor of
$$P:=g(g'_0,\ldots,g'_n)-f(g'_0,\ldots,g'_n)\,D\,y_{n+1},$$
where $g'_0,\ldots,g'_n\subset k[y_0,\ldots,y_n]$ define the inverse $\mathfrak{G}^{-1}$  to $\mathfrak{G}$
and $D$ is the target inversion factor.
In particular, one has $\deg(F)\leq \deg(f)\deg(\mathfrak{G}^{-1}) +\deg(D)+1$.
Moreover, the following conditions are equivalent:
\begin{enumerate}
\item[{\rm (a)}] $(P)=(F)$
\item[{\rm (b)}] $\gcd(g(\gg'), f(\gg')D)=1$
\item[{\rm (c)}] $\deg(F)= \deg(f)\deg(\mathfrak{G}^{-1}) +\deg(D)+1$.
\end{enumerate}
\end{Proposition}
\demo
Drawing on the above format of $\Psi$, consider the $1$-form corresponding to a Koszul syzygy as above
$$Q_i(\xx,\yy):= gy_i-fg_iy_{n+1},\; i\in\{0,\ldots,n\},$$
and take the corresponding syzygetic $\yy$-polynomial
$$P_i:=g(g'_0,\ldots,g'_n)y_i-f(g'_0,\ldots,g'_n)\,g_i(g'_0,\ldots,g'_n)\,y_{n+1}.$$
Note that $P_i$ is the numerator in the expression of $F$ as obtained in the proof of Proposition~\ref{pdegree} (i).

Clearly, then (a) through (c) are equivalent assertions.
\qed

\begin{Remark}\rm
One wonders what is a more precise choice of $f,g$ that guarantees the irreducibility of the form $P$ in
the above proposition.
Note that all Koszul-like syzygies of $J=(If,g)$ give rise to the same polynomial $P$, so there is not much
elbow room from this angle.
\end{Remark}

\section{The search for Rees equations}\label{third}

\subsection{The birational downgrading method}

For the results of this section, we recall a form of the so-called downgrading map in the context of birational maps.
Versions of this notion have been considered before in different contexts (\cite{BCS}, \cite{Trento}, \cite{syl2}).

Let $\xx,\yy$ be two sets of mutually independent variables over $k$ and of the same cardinality.
Given a bihomogenous polynomial $Q=Q(\xx,\yy)\in k[\xx,\yy]$, choose bihomogeneous
polynomials $Q_i(\xx,\yy), 0\leq i\leq n$, such that $Q=\sum_{i=0}^n x_i\,Q_i(\xx,\yy)$ -- called an
$\xx$-{\em framing} of $Q$.
In addition, fix a sequence of forms of the same degree $\mathfrak{H}:=\{\mathfrak {h}_0,\ldots, \mathfrak {h}_n\}\subset k[\yy]$.

The polynomial $\sum_{i=0}^n \mathfrak {h}_i\,Q_i(\xx,\yy)$ is called an
{\em  $\mathfrak{H}$-downgraded polynomial}
of $Q$.
We use the notation $D_{\mathfrak{H}}(Q)$ for an $\mathfrak{H}$-downgraded polynomial even though it is
not well-defined since the $\xx$-{\em framing} is only stable modulo the trivial (Koszul) relations
of $\xx$.
We will also allow for a harmless flat extension such as $k[\xx,\yy]\subset k[\xx,\yy,\zz]$, where $\zz$
is an additional set of variables.

This general notion will be applied to forms in $k[\yy,y_{n+1}]$ while $\mathfrak{H}\subset k[\yy]$ is the set of forms defining the
inverse of a Cremona map $\mathfrak{F}:\pp^n\dasharrow \pp^n$ -- in which case, we talk informally about a {\em birational downgrading}.
The common downgrading is typically the case where the Cremona map is the identity map.

As before, we stick to the notation $\mathcal{R}_R(J)\simeq R[\yy]/\mathcal{J}$
for the Rees algebra of an ideal $J\subset k[\xx]$ even if $\xx$ and $\yy$ have different
cardinalities.

\begin{Lemma}\label{downgrading_birational}
Let $\mathbf{g}=\{g_0,\ldots,g_n\}\subset R=k[\xx]$ be  forms of fixed degree defining a Cremona
map $\mathfrak{G}$ of $\pp^n$, not necessarily without a proper common divisor.
Let $g=g_{n+1}\in R$ stand for an additional form, of the same degree.
Write $\mathcal{J}\subset R[\yy,y_{n+1}]$ for the
presentation of the Rees algebra of the ideal $J=(\mathbf{g}, g)$ based on these generators.
Let $\mathfrak{H}:=\{g'_0,\ldots, g'_n\}\subset k[\yy]$ denote the
set of defining forms of the inverse map to $\mathfrak{G}$ and let $D_{\mathfrak{H}}$
denote the corresponding birational downgrading.
Then
$$Q=\sum_{r\geq 0} Q_r(\xx,\yy)y_{n+1}^r\in \mathcal{J}\, \Rightarrow\,D_{\mathfrak{H}}(Q)=\sum_{r\geq 0}
D_{\mathfrak{H}}(Q_r)y_{n+1}^r\in \mathcal{J}.$$
\end{Lemma}
\demo
Since the rational map $\pp^n\dasharrow \pp^{n+1}$ is birational onto its image
and $\mathfrak{H}$ (modulo the implicit equation) defines its inverse,
by a similar token as in (\ref{Rees_are_isomorphic}) one has an isomorphism
\begin{equation}\label{Rees_are_isomorphic2}
\mathcal{R}_R((\mathbf{g},g))\simeq \mathcal{R}_{S}((\mathfrak{H})),
\end{equation}
where $S=k[\yy,y_{n+1}]/(F)$ is the homogeneous coordinate ring of the image of $\mathfrak{F}$.
One proceeds as in the argument for a syzygetic polynomial, managing the respective
defining ideals.
However, instead of evaluating fully by $x_i\mapsto g'_i$, one only evaluates the  variables in a frame.
Since fully evaluating either $Q$ or its downgraded partner $D_{\mathfrak{H}}(Q)$ gives a form vanishing
on $\mathfrak{H}$ by the isomorphism (\ref{Rees_are_isomorphic2}), one has
$D_{\mathfrak{H}}(Q)\in \mathcal{J}$ as stated.
\qed

\bigskip

For the subsequent results we need an iterated version of the framing-downgrading gadget $D_{\mathfrak{H}}(Q)$.
Namely, one sets
\begin{equation}\label{framing_iteration}
D_{\mathfrak{H}}^0(Q)=Q, \quad D_{\mathfrak{H}}^{(\ell)}(Q):=D_{\mathfrak{H}}(D_{\mathfrak{H}}^{(\ell-1)}(Q)).
\end{equation}
We say that $D_{\mathfrak{H}}^{(\ell)}(Q)$ is {\em fully downgraded} when it eventually lands in $k[y_0,\ldots,y_{n+1}]$,
that is, when $\ell=\deg_{\xx}(Q)$.

We now apply to the original setup of the base ideal
$(If,g)\subset R$, where $I=(g_0,\ldots,g_n)$ is the base ideal of the Cremona map
$\mathfrak{G}$, and gcd$(f,g)=1$.
As before, let $g'_0,\ldots,g'_n$ have gcd $1$, defining the inverse map to $\mathfrak{G}$.
Accordingly, we take $\mathfrak{H}=\{g'_0,\ldots,g'_n\}$.
Note that, at least modulo Koszul relations, our previous syzygetic polynomials
are among the fully downgraded $D_{\mathfrak{H}}^{(\ell)}(Q)$, for $Q$ a syzygy of $J$ with nonzero last coordinate.

\begin{Proposition}\label{Rees_equations}
 The defining ideal of the Rees algebra of the ideal $J=(If,g)$
is a minimal prime of the ideal
$$\mathfrak{D}:=\left(\mathcal{I},\, \{D_{g'_0,\ldots,g'_n}^{(\ell)}(Q),\,\, 0\leq \ell\leq \deg_{\xx}(Q)\}\right),$$
where $\mathcal{I}$ stands for the defining ideal of $\mathcal{R}_R(I)$
and $Q\in k[\xx,\yy]$ runs through the biforms corresponding to the syzygies of
$J$ with nonzero last coordinate.
\end{Proposition}
\demo
By Lemma~\ref{downgrading_birational}, $\mathfrak{D}$ is contained in the presentation ideal of $\mathcal{R}_R(J)$
 which has codimension $n+1$ and is a prime ideal. Therefore, it suffices to show that
$\mathfrak{D}$ has codimension $n+1$ as well.
But $\mathcal{I}$ is a prime ideal of codimension $n$ and, moreover, is contained in the ideal $(\xx)k[\xx,\yy,y_{n+1}]$ because
$I$ is generated by algebraically independent elements over $k$. Since the fully downgraded elements
of $\mathfrak{D}$ belong to $k[\yy,y_{n+1}]$, this ideal has codimension at least one more.
\qed

\medskip

Concerning the problem of determining a set of generators of the presentation ideal of $\mathcal{R}_R(J)$, it is not enough to assume that $f$
is a general form in order that the (uniquely determined) fully downgraded be irreducible, as we have seen
in the non-zero-divisor case. But even when no Koszul relation is a minimal syzygy generator, taking $f$ general
may not help, as the following simple example indicates.

\begin{Example}\rm
Let $I=(x_0x_1,x_0x_2,x_1x_2)$ define the standard quadratic plane Cremona map.
Let $f=\lambda_0x_0+\lambda_1x_1+\lambda_2x_2$ be a general form (at least $\lambda_0\lambda_1\lambda_2\neq 0$).
Take $g=x_0^2x_1-x_2^3$, for example. The conductor $I:(g)$ is generated by $\{x_0,x_1\}$.
Accordingly, a set of minimally generating syzygies of $J=(If,g)$ consists of two linear syzygies coming from $I$
and two additional syzygies corresponding to $x_0,x_1$.
The syzygetic polynomial out of any of the two last syzygies has degree $5$
and has a so-called {\em extraneous} factor of degree $1$.
\end{Example}

\begin{Question}\rm Suppose that $f$ is a general form. Does a set of generators of the presentation ideal of $\mathcal{R}_R(J)$
consist of those of $\mathcal{I}$ plus the downgraded polynomials
$$\{D_{g'_0,\ldots,g'_n}^{(\ell)}(Q),\,\, 0\leq \ell\leq \deg_{\xx}(Q)\}$$
divided by the corresponding extraneous factors?
\end{Question}

The question lacks any precision since one would have to define ``extraneous factor''.
In any case, the ideal $\mathfrak{D}$ has a central place in this approach -- could be called
the {\em downgraded Rees ideal}.

\subsection{The method of the associated monoid parametrization}

In this subsection we will take a slightly different approach to get to the presentation ideal of
the Rees algebra of $J=(If,g)\subset k[\xx]$ defining a de Jonqui\`eres parametrization
$\mathfrak{F}:\pp^n\dasharrow \pp^{n+1}$, with underlying Cremona map $\mathfrak{G}:\pp^n\dasharrow \pp^{n}$.
As in the earlier notation, $I=(g_0,\ldots, g_n)\subset k[\xx]$, while the inverse map  $\mathfrak{G}^{-1}$
is defined by certain forms $g'_0,\ldots,g'_n\in k[\yy]=k[y_0,\ldots,y_n]$.

By Proposition~\ref{pdegree} (i), $F$ is an $y_{n+1}$-monoid, say, $F=F_{\delta}-y_{n+1}F_{\delta-1}$, where $\delta=\deg(F)$, and
$F_{\delta},F_{\delta-1}\in k[\yy]$ are forms of degrees $\delta,\delta-1$, respectively, such that $\gcd(F_{\delta},F_{\delta-1})=1$.

Set $h_{\delta}:= F_{\delta}(\xx)$ and $h_{\delta-1}:= F_{\delta-1}(\xx)$, so $h_{\delta}, h_{\delta-1}\in k[\xx]$ are forms
of degrees $\delta,\delta-1$, respectively.
Consider the standard monoid parametrization of ${\rm Im}(\mathfrak{F})$ defined by $h_{\delta}, h_{\delta-1}$, namely:

\begin{equation}\label{induced_monoidal}
\mathfrak{M}:=(h_{\delta-1}x_0:\cdots: h_{\delta-1}x_n:-h_{\delta}):\pp^n\dasharrow \pp^{n+1}.
\end{equation}

Write $K:=(h_{\delta-1}x_0,\ldots, h_{\delta-1}x_n,h_{\delta})\subset k[\xx]$ for the base ideal
of $\mathfrak{M}$.

\smallskip

Next is the main result of this part.

\begin{Theorem}\label{From_monoid}
With the above notation, one has:

\smallskip

\noindent {\rm (a)} $\mathfrak{F}$ and $\mathfrak{M}$ have the same implicit equation.

\smallskip

\noindent {\rm (b)} $\mathfrak{F}=\mathfrak{G}\circ \mathfrak{M}$.

\smallskip

\noindent {\rm (c)} Let
\begin{equation*}\label{presentations}
\mathcal{R}(J)\simeq k[\xx,\yy,y_{n+1}]/\mathcal{I}_{\mathfrak{F}}\quad {\rm and} \quad \mathcal{R}(K)\simeq
k[\xx,\yy,y_{n+1}]/\mathcal{I}_{\mathfrak{M}}
\end{equation*}
be presentations of the two Rees algebras based on the given generators.
Then
\begin{equation*}\mathcal{I}_{\mathfrak{F}}=\mathcal{I}_{\mathfrak{M}}(\mathfrak{G}):C^{\,\infty}\quad {\rm and}\quad \mathcal{I}_{\mathfrak{M}}=\mathcal{I}_{\mathfrak{F}}(\mathfrak{G^{-1}}):D^{\,\infty},
\end{equation*}
with
\begin{equation*}\mathcal{I}_{\mathfrak{M}}(\mathfrak{G}):=
\{\mathfrak{h}(g_0,\ldots,g_n; \yy, y_{n+1})\, |\, \mathfrak{h}(\xx; \yy, y_{n+1})\in \mathcal{I}_{\mathfrak{M}}\}
\end{equation*}
and
\begin{equation*}\mathcal{I}_{\mathfrak{F}}(\mathfrak{G^{-1}}):=
\{\mathfrak{h}(g'_0(\xx),\ldots,g'_n(\xx); \yy, y_{n+1})\, |\, \mathfrak{h}(\xx; \yy, y_{n+1})\in \mathcal{I}_{\mathfrak{F}}\}
\end{equation*}
where $C\in k[\xx]$ and $D\in k[\yy]$ are, respectively, the source inversion factor and the target inversion
factor of $\mathfrak{G}$.
\end{Theorem}
\demo (a) This follows immediately from Corollary~\ref{Benitez_Dandrea} and the definition of
the forms $h_{\delta}, h_{\delta-1}$.

\medskip

(b) This is pretty much tautological as the inverse of $\mathfrak{M}$ is induced (restriction) by the inverse of the identity
map of $\pp^n$ -- a special case of a de Jonqui\`eres parametrization, but otherwise very well-known (see, e.g.,
\cite{Piene_et_al}).
Then, obviously $\mathfrak{M}^{-1}\circ \mathfrak{F}=\mathfrak{G}$, as required.

\medskip

(c) We first show the inclusions $\mathcal{I}_{\mathfrak{F}}(\mathfrak{G^{-1}})\subset \mathcal{I}_{\mathfrak{M}}$ and
$\mathcal{I}_{\mathfrak{M}}(\mathfrak{G})\subset \mathcal{I}_{\mathfrak{F}}$.

For the first of these, let $\mathfrak{h}(\xx; \yy, y_{n+1})\in \mathcal{I}_{\mathfrak{F}}$ be a bihomogeneous
element.
By definition, one has $\mathfrak{h}(\xx; fg_0,\ldots,fg_n,g)=0$, while we wish to show that
$\mathfrak{h}(\gg'(\xx);h_{a-1}x_0,\cdots,h_{a-1}x_n,h_{a})=0$.
For this, let $D$ denote the target inversion
factor of $\mathfrak{G}$.
Recall that $h_{\delta}=\frac{g(\gg'(\xx))}{\deg(f)}$ and $h_{\delta-1}=\frac{f(\gg'(\xx))\,D(\xx)}{\deg(f)}$,
where $\deg(f)=\gcd(g(\gg'(\xx)), f(\gg'(\xx))\,D(\xx))$.
Since $\mathfrak{h}$ is bihomogeneous, we can pull out a power of $\deg(f)$ as a factor, hence the assertion is
equivalent to showing that
\begin{equation}\label{eqrees2}
\mathfrak{h}(\gg'(\xx);f(\gg'(\xx))D(\xx)x_0,\cdots,f(\gg'(\xx))D(\xx)x_n,g(\gg'(\xx)))=0.
\end{equation}
By definition, $D(\xx)x_i=g_i(\gg'(\xx)),\;\forall i$.  Therefore,  (\ref{eqrees2}) is equivalent to the vanishing of
$$\mathfrak{h}(\gg'(\xx);f(\gg'(\xx))g_0(\gg'(x)),\cdots,f(\gg'(\xx))g_n(\gg'(x)),g(\gg'(\xx)))=
\mathfrak{h}(\xx;fg_0,\cdots,fg_n,g)(\gg'(\xx)).$$
The rightmost polynomial is the resulting of evaluating the null polynomial, so itself is null.

\medskip

To argue for the second inclusion above, let likewise $\mathfrak{h}(\xx; \yy, y_{n+1})\in \mathcal{I}_{\mathfrak{M}}$
be a bihomogeneous element. By definition, $\mathfrak{h}(\xx; h_{\delta-1}\xx, h_{\delta})=0$ whereas one wishes to
show that
$$H:=\mathfrak{h}(\gg(\xx);f(\xx)\gg(\xx),g(\xx))=0.$$
For this, we first prove that substituting $\gg'(\xx)$ for $\xx$ in $H$ gives zero; namely, by a similar token
as above, using the characteristic property of the target factor $D$, there are suitable integers $s,r$ such that:
 \begin{eqnarray*}
H(\gg'(\xx))&=&\mathfrak{h}(\gg(\gg'(\xx));f(\gg'(\xx))g_0(\gg'(\xx)),\ldots, f(\gg'(\xx))g_n(\gg'(\xx)),g(\gg'(\xx)))\\
&=&\deg(f)^rD^s\mathfrak{h}\left(\xx;\frac{f(\gg'(\xx))}{\deg(f)}\gg(\gg'(\xx)),\frac{g(\gg'(\xx))}{\deg(f)}\right)\\
&=& \deg(f)^rD^s\mathfrak{h}(\xx;h_{\delta-1}\xx,h_{\delta})=0.
\end{eqnarray*}
Consider now the source inversion factor $C\in k[\xx]$, whose characteristic property is that $g'_i(\gg)=Cx_i,\;\forall i$.
Then, for a suitable exponent $t$, one has
$$C^{\:t}\,H= C^{\:t}\,\mathfrak{h}(\gg(\xx);f(\xx)\gg(\xx),g(\xx))=\left(\mathfrak{h}(\gg(\xx);f(\xx)\gg(\xx),g(\xx))\right)(\gg'(\gg(\xx)))$$
$$ =\left(\left(\mathfrak{h}(\gg(\xx);f(\xx)\gg(\xx),g(\xx))\right)(\gg'(\xx))\right)(\gg(\xx))=H(\gg'(\xx))(\gg(\xx)=0,$$
which proves the assertion.

\medskip

To complete the proof of the theorem, we  show the equality $\mathcal{I}_{\mathfrak{F}}=\mathcal{I}_{\mathfrak{M}}(\mathfrak{G}):C^{\,\infty}$,
the other equality being proved in the same fashion.
Since $\mathcal{I}_{\mathfrak{M}}(\mathfrak{G})\subset \mathcal{I}_{\mathfrak{F}}$ and the ideal $\mathcal{I}_{\mathfrak{F}}$
is prime, the inclusion $\mathcal{I}_{\mathfrak{F}}\supset\mathcal{I}_{\mathfrak{M}}(\mathfrak{G}):C^{\,\infty}$ is clear.
Conversely, let $\mathfrak{h}(\xx;\yy)\in \mathcal{I}_{\mathfrak{F}}$. By what we have proved above,
$\mathfrak{h}(\gg'(\xx);\yy)\in \mathcal{I}_{\mathfrak{M}}$ and hence $\mathfrak{h}(\gg'(\gg(\xx));\yy)\in
\mathcal{I}_{\mathfrak{M}}(\mathfrak{G})$.
Again,  $\gg'(\gg(\xx))=C\xx$ and $\mathfrak{h}(\xx;\yy)$ is bihomogeneous.
 Therefore, for a suitable exponent $u$, one has $C^{\:u}\mathfrak{h}(\xx;\yy)=\mathfrak{h}(\gg'(\gg(\xx));\yy)\in
 \mathcal{I}_{\mathfrak{M}}(\mathfrak{G})$, which says that $\mathfrak{h}(\xx;\yy)\in \mathcal{I}_{\mathfrak{M}}(\mathfrak{G}):C^{\,\infty}$.
This proves the other inclusion.
 \qed

\begin{Remark}\rm
The  result of Proposition~\ref{Rees_equations} and the one of Theorem~\ref{From_monoid} (c) give different
approaches to describe the presentation ideal $\mathcal{I}_{\mathfrak{F}}$ in an explicit way.
The first has the advantage of stressing a mechanical way to get the downgraded Rees ideal $\mathcal{D}$; unfortunately,
the final step may depend on the knowledge of the primary decomposition of $\mathcal{D}$.
The second has the advantage of starting with the simpler ideal $\mathcal{I}_{\mathfrak{M}}$, but is dependent
on knowing the implicit equation beforehand and a source inversion factor (the latter being equivalent, in practice,
to be able to get an inverse map explicitly).
This ideal has also been described in \cite[Theorem 3.1]{CorDandr} in a sort of ``reverse'' downgrading process
starting with the equation $F$.

It might be appropriate comparing the two procedures for the computational as well the theoretical purpose.
\end{Remark}

\medskip

\noindent Departamento de Matem\'atica, CCEN\\
Universidade Federal de Pernambuco\\
Cidade Universit\'aria, 50740-560 Recife, PE, Brazil.\\
aron@dmat.ufpe.br, hamid@dmat.ufpe.br

\end{document}